\documentclass[english]{amsart}

\usepackage{amsmath,amssymb,enumerate,bbm}

\usepackage[T1]{fontenc}
\usepackage[utf8]{inputenc}
\usepackage[all]{xy}

\usepackage{babel}
\usepackage{amstext}
\usepackage{amsmath}
\usepackage{amsfonts}
\usepackage{latexsym}
\usepackage{ifthen}

\usepackage{xypic}
\newtheorem{theo}{Theorem}[section]

\newtheorem{prop}[theo]{Proposition}
\newtheorem{lemm}[theo]{Lemma}
\newtheorem{coro}[theo]{Corollary}
\newtheorem{rema}[theo]{Remark}
\newtheorem{Defi}[theo]{Definition}
\newtheorem{ex}[theo]{Example}
\newtheorem{conj}[theo]{Conjecture}
\newtheorem{question}[theo]{Question}

\newcommand{\cqfd}
{%
\mbox{}%
\nolinebreak%
\hfill%
\rule{2mm}{2mm}%
\medbreak%
\par%
}
\newfont{\gothic}{eufb10}

\begin{document}
\title{Approximately rationally or elliptically connected varieties}
\author{Claire Voisin}
 \address{Institut de Math\'{e}matiques de Jussieu,
TGA Case 247, 4 place Jussieu, 75005 Paris, France}
\email{voisin@math.jussieu.fr} \maketitle \setcounter{section}{-1}
\begin{flushright} {\it To Slava Shokurov, on his 60th birthday}
\end{flushright}
\begin{abstract} We discuss a possible approach to the study of the vanishing of the Kobayashi pseudometric of a
projective variety $X$,  using
chains of rational or elliptic curves
contained in an arbitrarily small neighborhood of $X$ in projective space for the Euclidean topology.
\end{abstract}
\section{Introduction}
\setcounter{equation}{0}
In the paper \cite{lang}, Lang made some  conjectures concerning
entire curves in complex projective varieties $X$. He conjectured
for example that the Zariski closure of the locus in $X$ swept-out by entire curves is equal to the locus
swept-out by images of abelian varieties under  non constant rational maps $\phi:A\dashrightarrow X$.
When $X$ is a very general  quintic threefold in $\mathbb{P}^4$, this
has been shown  to be incompatible with Clemens conjecture \cite{clemensicm} by the following arguments : (i) $X$ contains countably
many rational curves, and
they are  Zariski   dense in $X$. (ii) On the other hand, $X$ is not swept out
by  images of non constant generically finite rational maps $\phi:A\dashrightarrow X$ with ${\rm dim}\,A\geq 2$ (cf. \cite{voisinscuola}).
(iii) Finally, if $X$ was swept-out by elliptic curves, this would contradict Clemens conjecture
on the discreteness of rational curves
in $X$ (see \cite[Lecture 22]{clemensinbook} or \cite{voisinscuola}).

The goal of this note is to discuss and illustrate by examples several  possible  notions of {\it approximate rational connectedness} or {\it approximate elliptic  connectedness} concerning complex projective manifolds. The general hope would be  that  approximately elliptically  connected
varieties are exactly varieties with trivial Kobayashi pseudodistance (cf. \cite{koba}). The main idea is that instead of looking at
rational or elliptic curves (or abelian varieties) sitting in $X$, we should study rational or elliptic
curves contained in arbitrarily small neighborhood of $X$ in projective space (for the Euclidean topology).

We assume  $X$ is embedded in some projective space  $\mathbb{P}^N$. We start with  the following na\"{\i}ve
definition:
\begin{Defi} \label{defi1} $X$ is said to be  approximately rationally connected in $\mathbb{P}^N$ {\bf in the na\"{\i}ve
sense}  if
for any neighborhood (for the Euclidean topology)  $U\subset \mathbb{P}^N$ of $X$, any two points of $X$ are
contained in a rational curve $C\subset U$.
\end{Defi}
\begin{rema} {\rm An equivalent definition is that any two points of $X$ can be joined
in an arbitrarily small neighborhood of $X$ by a chain of rational curves,
since such chains can be made irreducible by a small deformation in $U$, because $U$ has positive tangent bundle.}
\end{rema}

The reason why this definition cannot be interesting, from the point of view of
the study of the Kobayashi pseudodistance of $X$ is the following fact.
\begin{lemm} \label{exa10mai} i) Let $Y$ be a connected projective variety. Then $X:=Y\times \mathbb{P}^1$ is
approximately rationally connected in the na\"{\i}ve
sense in any projective embedding $X\subset \mathbb{P}^N$.

ii) More generally,  any connected variety $X$ such that the union of rational curves
contained in $X$ is dense in $X$ for the Euclidean topology is approximately rationally connected in the na\"{\i}ve
sense in any projective embedding.

iii) Assume $X\subset \mathbb{P}^N$ has the property that for any neighborhood $U$ of $X$, and for any point
$x\in X$, there is a rational curve $C\subset U$ passing through $x$. Then
$X$  is approximately rationally connected in the na\"{\i}ve
sense.
\end{lemm}
{\bf Proof.} i) Indeed, let $U$ be a neighborhood of $X$ in $ \mathbb{P}^N$.
Then for any automorphism $g$ of $\mathbb{P}^N$ sufficiently close to the identity,
and any curve $C_y=y\times \mathbb{P}^1\subset Y\times \mathbb{P}^1=X$, the curve $g(C_y)$ is contained in $U$.
It immediately follows that for any $x\in X$, the set of points $z$ in
$X$ such that there  is a rational curve $C\subset U$ passing through
$x$ and $z$ contains an open neighborhood   of the curve $C_{pr_1(x)}\setminus\{x\} $ in $ X\setminus\{x\}$.
Applying this argument to any $z\not=x$ in this neighborhood, we find that the set of points
$x'\in X$ such that there  is a chain of two rational curves $C_1,\,C_2 \subset U$ passing through
$x$ and $x'$ contains an open neighborhood   of $x $ in $ X$.
As $X$ is connected and compact, this easily implies
that  any two points of $X$ can be connected by a chain of rational curves contained in $U$.

ii) Let $U$ be a neighborhood of $X$ in $ \mathbb{P}^N$. For some $\epsilon>0$, $U$ contains
$U_\epsilon(X)=\{y\in \mathbb{P}^N,\,d(y,X)<\epsilon\}$. For any
point $x\in X$, there is by assumption a
 rational curve $C\subset X$ such that $d(x,C)<\epsilon$. Applying an automorphism $g$
 of $ \mathbb{P}^N$ such that $d(g,Id)<\epsilon$, we can thus find a curve $g(C)$ contained in $U$ and passing through
 $x$.
We then conclude using iii).

iii)  For any point
$x\in X$, there is a  rational  curve $C_x\subset U$ passing through $x$. Applying to them
 automorphisms  of $\mathbb{P}^N$  close to the identity and fixing a point $y\in C_x$, $y\not=x$,
we conclude as in i) that there is a  neighborhood $V_x$ of $x$ in  $X$ such that any point $y\in V_x$
is connected to $x$ by a
chain of two rational curves $C_{x}\cup g(C_{x})$ contained in $U$. By compactness and connectedness of $X$, we conclude
that any two points of $X$ can be joined by a chain of rational curves contained in $U$.
\cqfd

 \begin{rema} {\rm The statement in i) of the above lemma shows that the Kobayashi pseudodistance $d_{X,K}$
 of a subvariety
 $X\subset \mathbb{P}^N$ may  be different from
 the limit over the open sets $U\subset \mathbb{P}^N$ of the restrictions
 ${d_{U,K}}_{\mid X}$. Indeed,  in the above notation, if one chooses $Y$ to be Kobayashi hyperbolic, then the Kobayashi pseudodistance of $X=Y\times \mathbb{P}^1$ is non zero, while the
 restrictions ${d_{U,K}}_{\mid X}$ are all $0$.}

 \end{rema}
The main defect of  Definition \ref{defi1} is the fact that it is not stable under surjective
morphisms, that is, if $\phi:X\rightarrow Y$ is surjective and
$X$ is approximately rationally connected in the na\"{\i}ve sense, $Y$ needs not satisfy the same property. This follows indeed from
  Lemmas \ref{exa10mai}, i), and \ref{lekoba}.
  We could try to correct the definition by asking that not only $X$ but also
  all varieties $Y$ such that there is a surjective morphism
  from $X$ to $Y$, are approximately rationally connected in the na\"{\i}ve sense (say in any projective
  embedding).
  However the following example shows that this is not strong enough:
  \begin{ex}\label{ex23mai}{\rm
 Consider the case where  $X=(C\times S)/\iota$,
where $C$ is a curve of genus $\geq 2$ with hyperelliptic involution $i$, $S$ is a $K3$
surface which is the universal cover $S\rightarrow T$ of an Enriques surface,
the involution $\iota$ acting on
$C\times E$ acts as the hyperelliptic involution on $C$ and  as the involution over $T$  on  $S$.
 This involution $\iota$ has no fixed points. By Lemma \ref{exa10mai}, ii) $X$ is approximately rationally connected in any
 projective embedding, because rational curves are topologically dense in the fibers of $X\rightarrow
 \mathbb{P} ^1$. Consider any surjective morphism $X\rightarrow Y$, where $Y$ is normal.
  We claim that $Y$ is approximately rationally connected in the na\"{\i}ve sense in any projective embedding. Indeed,  if ${\rm dim}\,Y=1$, one has
 $h^{1,0}(Y)=0$, so $Y=\mathbb{P} ^1$. If ${\rm dim}\,Y=2$,  as $Y$ is dominated by $C\times S$, either it is  dominated
 by the $K3$ surface $S$,  or for each $c\in C$ the morphism from $c\times S$ to $Y$ has for image a curve $D$, and then, $Y$  is rationally
 dominated by a product $C\times \mathbb{P}^1$, since for any dominating rational map from a $K3$ surface to a smooth curve $D$,
 one has $D\cong \mathbb{P}^1$. In both cases, it  is approximately rationally connected in the na\"{\i}ve sense in any projective embedding, using
 Lemma \ref{exa10mai}. The case where ${\rm dim}\,Y=3$  is done similarly.}
 \end{ex}

 There are several ways to correct the na\"{\i}ve definition and we will propose two of them.
 \begin{Defi}\label{defi3} $X$ is  {\bf strongly} approximately rationally  connected if for any
  embedding $j$ of $X$ in a product $P$ of two projective spaces, $j(X)$
  is approximately rationally  connected
inside ${P}$  in the na\"{\i}ve sense (cf. Definition  \ref{defi1}).

\end{Defi}

Let us now give another variant of Definition \ref{defi1}, which might be easier to relate
to  the vanishing of the Kobayashi
pseudometric.
For any smooth $X\subset \mathbb{P}^N$, the projectivized tangent bundle
$\mathbb{P}(T_X)$ is naturally contained in the projectivized tangent bundle $\mathbb{P}(T_{\mathbb{P}^N})$. Any
curve $C\subset \mathbb{P}^N$ has a tangent lift $\widetilde{C}$ to $\mathbb{P}(T_{\mathbb{P}^N})$.
Let $\widetilde{U}$ be a neighborhood of $\mathbb{P}(T_X)$ in $\mathbb{P}(T_{\mathbb{P}^N})$. We will say that a curve
$C\subset \mathbb{P}^N$ is $\widetilde{U}$-close to $X$ if $\widetilde{C}\subset \widetilde{U}$. Hence, not only $C$ is close to $X$, but its tangent space at any point is close to $T_X$.
The following definition takes into account the cohomology classes of the curves considered. Here we
use the fact that if $U$ is a tubular neighborhood of $X$ in $\mathbb{P}^N$, $U$ and $X$ have the same homology. Let
 now $\widetilde{U}$ be a neighborhood of
$\mathbb{P}(T_X)$ in $\mathbb{P}(T_{\mathbb{P}^N})$.
Note that if   for any $U,\,\widetilde{U}$, and for any point $x$ of $X$,  there passes a curve $C_x\subset {U}$ passing through
$x$, which is
$\widetilde{U}$-close to $X$, $C_x$ can be chosen to vary locally continuously with $x$, hence to have a cohomology class
$[C]\in H_2(U,\mathbb{Z})=H_2(X,\mathbb{Z})$ locally independent of $x$. By considering chains, and by smoothing them, we conclude using arguments similar to the proof of Lemma \ref{exa10mai}, that if
$X$ is connected, we can assume that the class of the covering curves $C_x$  is in fact  independent of $x$.

\begin{Defi}\label{defi5} A connected variety $X\subset \mathbb{P}^N$ is  {\bf cohomologically} approximately rationally  connected if  for any $U,\,\widetilde{U}$ as above, through any point $x$ of $X$ there passes a rational curve $C_x$ contained in $U$ and $\widetilde{U}$-close to
$X$, of class $[C]$ independent of $x$, and furthermore
the convex cone generated by the $(n-1,n-1)$-components of
classes $[C_i]\in H_2(U,\mathbb{Z})=H_2(X,\mathbb{Z})=H^{2n-2}(X,\mathbb{Z})$
 of such covering curves $C_{i,x}$  contains a strongly positive class.

\end{Defi}
Here we say  that a   class  of type $(n-1,n-1)$ on an $n$-dimensional variety $X$  is strongly positive
if it  has positive intersection with  pseudoeffective  $(1,1)$-classes (represented by
weakly positive currents of type $(1,1)$). When the class belongs to the space $N_1(X)$ generated by
curve classes, this is equivalent to being  in the interior of the convex cone generated by classes of moving curves (cf. \cite{bdpp}).
\begin{rema}{\rm Since the class of a curve $C\subset U$ is the class of the
current of integration over $C$, this is the class in $U$ of a closed current of type $(N-1,N-1)$, $N={\rm dim}\,U$.
Using a differentiable retraction $\pi:U\rightarrow X$, we also have
 the current of integration over $\pi(C)$, whose class is the cohomology class
$[C]$
above. This last class is not in general of type $(n-1,n-1)$ (see examples in section \ref{sec1}). However, when $\widetilde{U}$ is small,
 it is  close to be of type $(n-1,n-1)$, as  $C$ is  $\widetilde{U}$- close to
$X$.
}
\end{rema}
The cohomological condition in Definition \ref{defi5}
 addresses the weakness of Definition \ref{defi1}: indeed, the rational
curves $g(C_y)$ used in the proof of Lemma \ref{exa10mai}, i) are in the same class as
the fibers of $pr_1$, and this class is not strongly positive.

\begin{rema}{\rm If $H^2(X,\mathbb{Q})=\mathbb{Q}$, Definitions \ref{defi5} and \ref{defi1} are quite close. Indeed, in this case,
the cohomological condition in Definition \ref{defi5} is empty, and we thus just ask that for any
neighborhoods $U$ of $X$, $\widetilde{U}$ of $\mathbb{P}(T_{X})$, and any general point $x\in X$, there is a rational curve in $U$ passing through $x$ and $\widetilde{U}$-close to $X$.
The last condition is satisfied by the examples of Lemma \ref{exa10mai}, i) (but they do not satisfy $H^2(X,\mathbb{Q})=\mathbb{Q}$).}
\end{rema}
We believe that these notions should be related to the triviality of
the Kobayashi pseudodistance   (see \cite{koba}) of $X$, although it is quite hard
to establish precise relations. This is due  to the notorious
difficulty to localize
Ahlfors currents or Brody curves (see \cite{duval}, \cite{duvalens}, \cite{paun} for important progresses
on this subject). The motivation for introducing these geometric definitions is the lack of progress on the understanding of complex varieties with vanishing Kobayashi pseudodistance
(by contrast to the recent progresses made on the Green-Griffiths conjecture, eg  for high degree hypersurfaces in projective space, see \cite{DMR}).

However, the following easy lemma shows that
approximate rational connectedness in either of the above senses is too restrictive topologically:
\begin{lemm} \label{letopo} Abelian varieties are not approximately rationally connected (in the na\"{\i}ve sense). More precisely, if
$A\subset \mathbb{P}^N$ is an abelian variety, and $U$ is a tubular neighbourhood of $A$,
$U$ does not contain any rational curve.
\end{lemm}
{\bf Proof.}   $A$ and $U$ have the same homotopy
type.
Hence, as $\pi_2(A)=0$, we also have $\pi_2(U)=0$. Thus a rational curve contained in $U$ should be homologous to $0$ in $U$,
hence in $\mathbb{P}^N$, which is absurd.
\cqfd
This lemma (and the fact that abelian varieties have trivial Kobayashi
pseudodistance) is the motivation for the following variant of Definitions \ref{defi3}, \ref{defi5},
\ref{defi1}:
\begin{Defi} \label{defi6} i) $X\subset \mathbb{P}^N$ is said to be approximately elliptically connected in the na\"{\i}ve sense if for any
neighborhood $U$ of $X$ in $\mathbb{P}^N$,
any two points $x,y\in X$ can be joined by a chain of elliptic curves in $U$.

ii) $X\subset \mathbb{P}^N$ is said to be strongly
approximately elliptically connected if Definition \ref{defi3}  holds with rational curves replaced by chains of elliptic curves.

iii) $X\subset \mathbb{P}^N$ is said to be
 cohomologically  approximately elliptically connected if Definition  \ref{defi5}
holds with rational curves replaced by  elliptic curves.

\end{Defi}

\begin{rema} {\rm  If $X$ is connected and  approximately  rationally (resp. elliptically) connected in the  cohomological
sense, then it is also approximately   rationally (resp. elliptically) connected in the na\"{\i}ve sense, as follows from
Lemma \ref{exa10mai}, iii) (or its obvious extension to the elliptic case).}
\end{rema}

Properties i), ii) and iii) are    satisfied by (very general) Calabi-Yau varieties obtained as the double cover of
projective space $\mathbb{P}^n$ ramified along a degree $2n+2$ hypersurface as they are
covered (in infinitely many different ways)
by families of elliptic curves (cf. \cite{voisincurrent}).

We will give however in section \ref{sec1} (see Theorem \ref{exampletheo})
examples of varieties containing only finitely many rational curves, but which are approximately
 rationally connected in the cohomological sense. Similarly, abelian varieties are approximately  cohomologically
elliptically connected
(see Theorem \ref{propexample}), while the general ones do not contain any elliptic curve. From
this one can deduce that  Fano varieties of lines of very general cubic threefolds
 satisfy this property, as they are covered in infinitely many different ways by a $2$-dimensional family of surfaces
birationally equivalent to abelian surfaces
(cf. \cite{voisinfano}).

Our hope  is that approximately  elliptically  connected varieties in one of the strengthened senses described above are
the same as the ``special varieties'' invented by Campana \cite{cafou} (which are also conjectured
to be the complex projective varieties with trivial Kobayashi pseudodistance).
Notice that Demailly in \cite{demailly} gives a description of the Kobayashi pseudometric
of $X$
involving algebraic curves in $X$, together with their intrinsic hyperbolic metric.
 This says that if $X$ has a trivial Kobayashi pseudodistance, there are many
 algebraic curves
 in $X$ for which the intrinsic hyperbolic metric is small compared to the
 metric obtained by restricting a given metric on $X$. In particular this compares
 the genus of these curves to their degrees, but
 this does not say anything on  the genus alone.

 To start with, we have the following easy lemma.
\begin{lemm} \label{lekoba} If a projective variety $X$ is Kobayashi hyperbolic, it is not approximately elliptically  connected
(in the na\"{\i}ve sense) in any projective embedding.
\end{lemm}
{\bf Proof.} Indeed, if there is an elliptic curve $E_n$ in any neighborhood $U_{\epsilon_n}(X)$ of
$X$ in $\mathbb{P}^N$, with $\lim_{n\rightarrow\infty}\epsilon_n=0$, we can choose
for each $n$  a holomorphic map
$f_n:\Delta\rightarrow E_n\rightarrow \mathbb{P}^N$  from the unit disk to $E_n$, such that
$|f_n'(0)|=n$, where the modulus of the derivative is computed with respect to the ambient metric.
 By Brody's lemma \cite{brody}, there is an entire curve in $\mathbb{P}^N$  obtained as the limit of
 a subsequence of the $f_n$'s conveniently reparametrized. This entire curve is contained in
 $\cap_nU_{\epsilon_n}(X)=X$ and $X$ is not Kobayashi hyperbolic.

\cqfd
We will also prove in
section \ref{sec2}  the following property:
\begin{prop} (see Proposition \ref{theosurj})
\label{propsurj} If $X$ is strongly approximatively rationally  (resp. elliptically) connected  and
$\phi: X\rightarrow Y$ is a surjective morphism, then $Y$ is approximately rationally  (resp. elliptically) connected
in the na\"{\i}ve sense.  In particular, $Y$ is not Kobayashi hyperbolic.
\end{prop}

We do not know whether this result holds for  cohomological approximate elliptic or rational connectedness, but we know by
Lemmas \ref{exa10mai} and \ref{lekoba} that it does not hold for na\"{\i}ve approximate elliptic or rational connectedness.

Next, as  the definitions are stable under \'{e}tale  covers (see Lemma \ref{letale}), one crucial point
needed in order
to  make the class of approximately elliptically connected varieties close to Campana's special manifolds would be the following:
\begin{conj} A variety of general type is not approximately elliptically  connected in the na\"{\i}ve sense.
\end{conj}
As we do not even know that elliptic or rational curves are not topologically dense in
a variety of general type  (a weak version of Green-Griffiths-Lang conjecture), this question seems to
be out of reach at the moment. For example, we  know that general hypersurfaces in $\mathbb{P}^{n}$
of degree $\geq 2n-2$ do not contain any rational curve for $n\geq 4$ (see \cite{voisinrat}, giving  an optimal bound which is  slightly better
than  \cite{clemens})
and that the only rational curves contained in general hypersurfaces in $\mathbb{P}^{n}$
of degree $ 2n-3$, for $n\geq 6$, are lines (cf. \cite{pacienza}),
but general hypersurfaces in $\mathbb{P}^{n}$ of degree $n+2\leq d\leq 2n-4$ are not known to
carry finitely many families of rational curves and not even known not to contain a dense set covered
by rational curves.

In the other direction, we do not know if rational curves in Calabi-Yau hypersurfaces
are topologically dense, except in dimension $2$, that is for $K3$ surfaces,  for a Baire second category subset of the
moduli space by \cite{chenlewis}. One question implicitly raised in the present paper is whether it is easier to study rational (or elliptic) curves contained in small neighbourhoods of such hypersurfaces.

\vspace{0.5cm}

{\bf Thanks.} {\it  I thank Fr\'{e}d\'{e}ric Campana, Tommaso de Fernex, J\'anos Koll\'ar, Mihai P\u{a}un and Jason Starr
for useful discusions related to this subject.
This paper was written during my stay at Isaac Newton Institute, where I benefitted from ideal working conditions. I warmly thank  Peter Newstead, Leticia Brambila-Paz,
 Oscar Garc\'{\i}a-Prada and  Richard Thomas
for their invitation to participate in  the semester on Moduli Spaces. }

\section{\label{sec1} Some examples}
 Let us give two examples of classes of  varieties which do not contain many rational (resp. elliptic curves)
 but are approximately rationally (resp. elliptically) connected in the cohomological sense.

 \begin{theo}  \label{propexample} Abelian varieties are approximately elliptically connected
 (for any projective embedding) in the  cohomological sense.
 \end{theo}
 {\bf Proof.} Let $A\subset \mathbb{P}^N$ and let $\widetilde{U}\subset \mathbb{P}(T_{\mathbb{P}^N})$ be an Euclidean neighborhood of
 $\mathbb{P}(T_A)$.
  For a   small deformation $A_\epsilon$ of $A$ in $\mathbb{P}^N$,
  $\mathbb{P}(T_{A_\epsilon})\subset \mathbb{P}(T_{\mathbb{P}^N})$ is a small deformation of
  $\mathbb{P}(T_A)$, hence will be contained in $\widetilde{U}$ when the deformation is small enough. It follows
  that for any curve $C\subset A_\epsilon$, its tangent lift $\widetilde{C}$
  is contained in $\widetilde{U}$.

We use now  the well-known fact that abelian varieties isogenous to a product $E_1\times\ldots\times E_n$, where each $E_i$ is an elliptic curve, are dense
(for the Euclidean topology) in the moduli space of  $n$-dimensional polarized abelian varieties. On the other hand,
inside $E_1\times\ldots\times E_n$, the elliptic curves obtained as the images of $E_i$ under the natural morphisms
$\phi_i:E_i\rightarrow E_1\times\ldots\times E_n,\,x\mapsto (e_1,\ldots,e_{i-1},x,e_{i+1},\ldots, e_n)$, for given points $e_j\in E_j$ can be chosen to pass through any point. Of course, the same is true
 for any abelian variety isogenous to $E^n$.

Let now $A_\epsilon\subset \mathbb{P}^N$ be a sufficiently small deformation
of $A$ which is isogenous to $E_1\times\ldots\times E_n$. Then the elliptic curves
$\phi_i(E_i)$
contained in $A_\epsilon$ sweep-out $A_\epsilon$ and their tangent lift
is contained in $\widetilde{U}$. For any point $x\in A$, we can find an automorphism
of $\mathbb{P}^N$ which is close to the identity and such that $x\in g(A_\epsilon)$. Thus,
the curves $g(\phi(E))$ can be chosen to pass through any point of $A$ and to have their tangent lift
contained in $\widetilde{U}$.

Finally, when there is no nonzero morphism between $E_i$ and $E_j$ for $i\not= j$, the classes of the curves $\phi_i(E_i)$ generate
the space of Hodge classes $Hdg^{2n-2}(A_\epsilon)$.  In particular, a convex combination
of these classes contains the class $h_{A_\epsilon}^{n-1}$, where $h_{A_\epsilon}=c_1(\mathcal{O}_{A_\epsilon}(1))$.
Using the canonical  isomorphism $ H^{2n-2}(A_\epsilon,\mathbb{R})\cong  H^{2n-2}(A,\mathbb{R})$ we conclude
that a convex combination
of these classes, transported to $A$, contains the class $h_A^{n-1}$, where $h_A=c_1(\mathcal{O}_{A}(1))$. Taking the $(n-1,n-1)$-component, we conclude as well hat a convex combination
of the $(n-1,n-1)$-components of these classes  transported to $A$  contains the class $h_A^{n-1}$, which finishes the proof.

\cqfd

\begin{rema}{\rm The example of abelian varieties also illustrates why the notion of approximate
elliptic (or rational) connectedness might be easier to study than the
property of  being swept out by entire curves,
or of having arbitrarily small  neighborhoods in $\mathbb{P}^N$ swept-out by  entire curves. Indeed, there are of course
a lot of entire curves in abelian varieties. However, the elliptic curves $E$ exhibited above, contained
in a small deformation of a given abelian variety $A$ in projective space, are much more reasonable, since
their induced metric  is uniformly equivalent to their flat metric $k_E$ (normalized so that the volume is equal to the
 degree). This is because the flat metric of $A$ is equivalent
to the induced metric on $A$, which  easily implies that
there is a flat metric $h_{A_\epsilon}$ on $A_\epsilon$  equivalent
to the induced metric on $A_\epsilon$  with constants   depending only on $A$.
The restriction $h_{A_\epsilon\mid E}$ of this flat metric to $E$ is a flat metric on $E$ which is then
uniformly equivalent to the induced metric $h_{\mid E}$.
In other words we have $c h_{\mid E} \leq h_{A_\epsilon\mid E}\leq Ch_{\mid E}$, for some constants
$c,\,C$ depending only on $A$.
Integrating over $E$ the corresponding K\"ahler forms, we get
$c\,{\rm deg}\,E\leq \int_E\omega_{A_\epsilon\mid E}\leq C\,{\rm deg}\,E$,
which tells, since $\omega_{A_\epsilon\mid E}$ is the flat metric, that the normalized
metric $k_E$, which is equal to $\frac{deg\,E}{\int_E\omega_{A_\epsilon\mid E}}h_{A_\epsilon\mid E}$
satisfies
\begin{eqnarray}\label{equiv} \frac{c}{C} h_{\mid E} \leq k_E\leq \frac{C}{c}h_{\mid E}.
\end{eqnarray}

}
\end{rema}
One interesting question is the following:
\begin{question} \label{questionmetric}i)  Does  any elliptic curve close enough in the usual topology to an abelian
subvariety $A$ of $\mathbb{P}^N$ satisfy (\ref{equiv}) for some constants depending only on $A$?

 ii) Does the above question have an affirmative answer  for elliptic curves
$\widetilde{U}$-close to $A$, for a small neighborhood $\widetilde{U}$ of $\mathbb{P}(T_A)$ in $\mathbb{P}(T_{\mathbb{P}^N})$?
\end{question}
An affirmative answer to these questions would have the following consequence:
\begin{prop} \label{proputopie} Assume Question \ref{questionmetric}, i)
has an affirmative answer for a given abelian variety $A\subset \mathbb{P}^N$. Then
a subvariety $X\subset A$ which is of general type is not approximately elliptically connected in the na\"{\i}ve sense. If Question  \ref{questionmetric}, ii) has an affirmative answer for  $A\subset \mathbb{P}^N$, then
a subvariety $X\subset A$ which is of general type is not approximately elliptically connected in the cohomological sense.
\end{prop}
{\bf Proof.} Indeed, one knows by \cite{kawamata} that $X$ satisfies the Green-Griffiths conjecture, so that the union
of the entire curves contained in $X$
is not Zariski dense in $X$. On the other hand, assume that for any $x\in X$, there is
an elliptic curve $E_n\subset \mathbb{P}^N$ passing through $x$ such that $E_n\subset V_{\frac{1}{n}} (X)=\{y\in \mathbb{P}^N,\,d(y,X)<\frac{1}{n}\}$.
Then consider the flat uniformization $f_n:\mathbb{C}\rightarrow E_n$ such that
$f_n(0)=x_n$ and $f_n^*k_{E_n}$ is the standard metric (so $f_n$ is defined up to the action
of $U(1)$). Then
if  Question \ref{questionmetric}, i)  has a positive answer, as $E_n$ is close to
 $A$, the  derivatives $|f_n'|$ (computed with respect to the ambient metric on
$\mathbb{P}^N$) are bounded above and below in modulus, so that
we can extract a subsequence which converges uniformly on compact sets of $\mathbb{C}$ to a non constant
entire curve passing through $x$ and contained in $X$. As $x$ is arbitrary, this contradicts Kawamata's result.
Similarly, if question ii) has an affirmative answer, elliptic curves $\widetilde{U}$-close to
$A$ satisfy the above property for $\widetilde{U}$ small. This is then also true
for elliptic curves $\widetilde{U}$-close to $X$. Hence, by the above argument, one cannot have
an elliptic curve $\widetilde{U}$-close to $X$ passing through any point of $X$ for $\widetilde{U}$
arbitrarily small.
\cqfd
The second example we will consider is the example of elliptic surfaces with finitely many rational curves. More precisely, we consider a very general hypersurface $S\subset \mathbb{P}^1\times \mathbb{P}^2$ of bidegree $(2l,3)$ with $l\geq 2$.
\begin{theo} \label{exampletheo} i) $S$ contains finitely many rational curves, namely the singular fibers of the
elliptic fibration $f:=pr_{1\mid S}:S\rightarrow \mathbb{P}^1$.

ii) $S$ is approximately rationally connected (relative to the Segre embedding) in the  cohomological sense.
\end{theo}
{\bf Proof.} i) A smooth surface $\Sigma\subset \mathbb{P}^1\times \mathbb{P}^2$ of bidegree $(2,3)$
is a $K3$ surface, which contains only countably many rational curves. If $\Sigma$ is chosen to be very general,
$${\rm Pic}\,\Sigma=({\rm Pic}\,(
\mathbb{P}^1\times \mathbb{P}^2))_{\mid \Sigma}$$
hence no algebraic curve in $\Sigma$ has degree $1$ over $\mathbb{P}^1$.
Take such a $\Sigma$ and consider a very general morphism    $\phi:\mathbb{P}^1\rightarrow \mathbb{P}^1_0$ of degree
$l$. Then the surface $S=\Sigma\times_{\mathbb{P}^1_0}\mathbb{P}^1$ is of bidegree $(2l,3)$
in $\mathbb{P}^1\times \mathbb{P}^2$ and for each rational curve $C\subset \Sigma$ not contained in a fiber of $\Sigma\rightarrow\mathbb{P}^1_0$, the
curve $C\times_{\mathbb{P}^1_0}\mathbb{P}^1\subset S$ is irreducible of positive geometric genus.
Hence $S$ does not contain any other rational curves than those  contained in a fiber of $pr_1:S\rightarrow \mathbb{P}^1$.

ii) First of all, we apply the criterion for density
of the Noether-Lefschetz locus due to M. Green (cf. \cite[Proposition 5.20]{voisinbook}) to show that
arbitrarily small deformations
$S_\epsilon\subset \mathbb{P}^1\times \mathbb{P}^2$
of a general surface $S$ as above admit sections of $f_\epsilon:\,S_\epsilon\rightarrow \mathbb{P}^1$.

We recall the statement of the criterion in the form which we will  use here : consider the universal family
$\pi:\mathcal{S}\rightarrow B,\,\mathcal{S}\subset B\times \mathbb{P}^1\times \mathbb{P}^2$ of such smooth surfaces. Let $0\in B$ and let $S_0$ be the fiber $\pi^{-1}(0)$.
Then we have:
\begin{prop} \label{density}Assume there is a $\lambda\in H^1(S_0,\Omega_{S_0})$ such that the
map $\overline{\nabla}_\lambda:T_{B,0}\rightarrow H^2(S_0,\mathcal{O}_{S_0})$ is surjective.
Then for any open set $U\subset B$ (for the Euclidean topology),
the set of classes $\alpha\in H^2(S_0,\mathbb{Z})$ which become algebraic on some fiber $\mathcal{S}_t$ for some $t\in U$ contains the set of integral points in a non-empty open cone of $H^2(S_0,\mathbb{R})$.
\end{prop}
In this statement, the map $\overline{\nabla}_\lambda$ is the composition of the Kodaira-Spencer map
$T_{B,0}\rightarrow H^1(S_0,T_{S_0})$ and the cup-product/contraction map
$\lambda\lrcorner: H^1(S_0,T_{S_0})\rightarrow H^2(S_0,\mathcal{O}_{S_0})$.

Using the Griffiths description of the IVHS of hypersurfaces (see \cite[6.2.1]{voisinbook}), we find that
this map identifies to the  multiplication
$$\mu_{P_\lambda}: R_{2l,3}(S_0)\rightarrow R_{6l-2,6}(S_0)$$
by a certain polynomial
$P_\lambda\in H^0(S_0,\mathcal{O}_{S_0}(4l-2,3))$, where $R(S_0)$ is the Jacobian ring of the defining
equation of $S_0$.
 One checks explicitly that for generic
$S_0$ and generic $P_\lambda$, the map $\mu_{P_\lambda}$ is surjective.

Using Proposition \ref{density}, we conclude that for generic $S_0$ and
 for any small simply connected neighborhood $U$ of $0$ in
$B$, there is a non-empty open cone $C$ in $H^2(S_0,\mathbb{R})$ such that
any integral class $\alpha\in H^2(S_0,\mathbb{Z})\cap C$ becomes (by parallel
transport) algebraic on some fiber
$S_t$, for some $t\in U$.

Let now $F\in H^2(S_0,\mathbb{Z})$ be the class of a fiber of $f_0$. Then  its parallel transport
to $S_t$ is the class $F_t$ of a fiber of $f_t$ and the elliptic  fibration $f_t:S_t\rightarrow
\mathbb{P}^1$ admits a section if and only if there is an algebraic class
$\alpha\in H^2(S_t,\mathbb{Z})\cong H^2(S_0,\mathbb{Z})$ such that $<\alpha,F>=1$.
Note furthermore that the class $h:=c_1(pr_2^*\mathcal{O}_{\mathbb{P}^2}(1))$, which is of degree $3$ on the class $F$, remains algebraic
on any deformation $S_t$ of $S_0$ in $\mathbb{P}^1\times \mathbb{P}^2$.
Hence, $f_t$ has a section if and only if there is an algebraic class
$\alpha\in H^2(S_t,\mathbb{Z})\cong H^2(S_0,\mathbb{Z})$ such that $<\alpha,F>=1$ mod. $3$.
To conclude that the set of surfaces $S_t$ having a section
is dense, we then use  the following lemma (which is used implicitly in  \cite[Remark 1]{voisinsoule}).
\begin{lemm} \label{leF1dense} For any non-empty open cone $C\subset H^2(S_0,\mathbb{R})$, there are
 elements in $C\cap \{\alpha\in H^2(S_0,\mathbb{Z}),\,<\alpha,F>=1\,\, mod. \,\,3\}$.
\end{lemm}

The second and final step of the proof is the following lemma, due to Chen and Lewis
 \cite{chenlewis}:
\begin{lemm} \label{lexile}Let $f:S\rightarrow \mathbb{P}^1$ be an elliptic fibration, and $L$ a
line bundle on $S_t$, of degree $d\not=0$ on fibers. Assume the fibers of $f$ are irreducible and reduced and that
the monodromy group $\pi_1(\mathbb{P}^1_{reg},t_0)\rightarrow Aut \,H^1(S_{t_0},\mathbb{Z})$ is the full
symplectic group ${\rm SL}(2,\mathbb{Z})$. Then for any section $\sigma :\mathbb{P}^1\rightarrow S$ of $f$ such that
the class $d[\sigma]-c_1(L)$ is not of torsion, the curves
$C_n:=\sigma_n(\mathbb{P}^1)$ have the property that $\cup_nC_n$ is  dense for the Euclidean topology
in $S$.
\end{lemm}
Here $\sigma_n:=\mu_n\circ \sigma$, where $\mu_n:S\dashrightarrow S$ is the
self-rational map which to $x\in  S$ with $f(x)=u\in \mathbb{P}^1$ associates
the point $y$ of the fiber $S_u$ such that $(dn+1)x-nL_{\mid S_u}=y$ in ${\rm Pic}\,S_u$.

\cqfd

 The  proof of Theorem \ref{exampletheo} is now  concluded as follows:
 let $S_0$ be generic. Let $U$  be a neighborhood of $S_0$ in $ \mathbb{P}^5$ and  $\widetilde{U}$ be a neighborhood of
 $\mathbb{P}(T_{S_0})$ in $\mathbb{P}(T_{\mathbb{P}^5})$.
 As $S_0$ is generic, the fibration $f_0:S_0\rightarrow \mathbb{P}^1$
 is a Lefschetz fibration with irreducible fibers and with monodromy
 equal to the full symplectic group. These properties remain true for
 a small deformation of $S_0$.

 By Lemma \ref{leF1dense}, there is a surface $S_t$ which is a small deformation of $S_0$
  in $\mathbb{P}^1\times \mathbb{P}^2\subset \mathbb{P}^5$
   such that $f_t:S_t\rightarrow \mathbb{P}^1$
  has a section $C$. In particular $S_t$ is contained in $U$ and  $\mathbb{P}(T_{S_t})$ is contained in $\widetilde{U}$. Thus any curve contained in
  $S_t$ is $\widetilde{U}$-close to $S_0$.
  By lemma \ref{lexile}, the union $\cup_nC_n$ is dense in $S_t$. Applying automorphisms of
  $\mathbb{P}^5$ close to the identity if needed, we conclude that for any point
  $x\in S_0$, there is an arbitrarily small deformation $S_t$ of $S_0$ in $\mathbb{P}^5$ containing a section $C_n$ passing through
  $x$. These curves are contained in $U$ and $\widetilde{U}$-close to $S_0$.

  To conclude, it only remains to check that a convex combination of the $(1,1)$-components of the classes of (a variant of) the curves $C_n$
  (transported to $S_0$) contains an ample class on $S_0$. For this we observe that the sum $2[C_n]+[C_{-2n-1}]\in H^2(S_t,\mathbb{Z})=H^2(S_0,\mathbb{Z})$ is a combination
  of the class $h$ and the class $F$. The coefficient in $h$ is obviously positive.
   This class may not be ample, but we observe that the class $F$ is in $S_0$ the class of a rational curve (namely a singular fiber).
  Instead of the curves $C_n$, we can thus consider  curves $C'_n$ obtained by smoothing in $\mathbb{P}^5$
  the union of $C_n$ and a covering of large degree of a singular fiber. The resulting curves can be chosen to stay $\widetilde{U}$-close to $S_0$
  and the sum $2[C'_n]+[C'_{-2n-1}]\in H^2(S_t,\mathbb{Z})=H^2(S_0,\mathbb{Z})$ is an ample class.

\cqfd

These two examples might give the feeling that
the  natural way to produce elliptic or rational curves in
a (arbitrarily) small neighborhood of a subvariety $X\subset \mathbb{P}^N$
is by studying elliptic or rational curves lying in some small deformation
$X_\epsilon$ of $X$ in $\mathbb{P}^N$. This is however not true at all, as shows
 the following example, obtained by mimicking the trick of \cite{voisininvent}: start with an abelian variety $A\cong \mathbb{C}^n/\Gamma$ admitting an endomorphism $\phi$
such that the corresponding endomorphism $\phi_\mathbb{Q}$ of $\Gamma_\mathbb{Q}=H_1(A,\mathbb{Q})$
has only eigenvalues of multiplicity $1$. Then it is immediate to see that
the pair $(A,\phi)$ is rigid. Furthermore, we can construct such $\phi$'s so that
$A$ does not contain any elliptic curve (one considers for example  simple abelian varieties with complex multiplication).

Starting from such pair $(A,\phi)$ we consider the projective variety $X$ obtained by blowing-up
successively $A\times A\times \mathbb{P}^1$ along $A\times x_0\times t_1,\,x_1\times A\times t_2,\, {\rm diag}\,A\times t_3,\,{\rm graph}\,\phi\times t_4$, for generic choices of $x_0,\,x_1$ and distinct points $t_1,\ldots, t_4\in \mathbb{P}^1$.
Choose a projective  embedding of $X$ in $\mathbb{P}^N$.
As $A$ does not contain elliptic or rational curves, the only rational or elliptic curves contained in $X$ are contained in the union $D$ of the
exceptional divisors
of the blow-ups or in proper transforms of the fibers of the map $p_{A\times A}\circ \tau: X\rightarrow A\times A$, where $\tau:X\rightarrow A\times A\times \mathbb{P}^1$ is the blow-up map.
 Furthermore, the deformations of $X$ preserve the exceptional divisors
hence are all of the same type as $X$, and as the pair $(A,\phi)$ is rigid, it follows that
elliptic or rational curves contained in a small deformation $X_\epsilon$ of $X$ in
$\mathbb{P}^N$ are close (for the usual topology) to either
a curve  contained
in $D_\epsilon$, or to a fiber of the map $p_{A\times A}\circ \tau: X\rightarrow A\times A$.
Hence, for a general point $x$ of $X$, elliptic curves passing through $x$ and contained in a small deformation
$X_\epsilon$ of $X$ have for  homology  class  a multiple of the class of a fiber of  $p_{A\times A}\circ \tau$.  As this class
  is not strongly positive, we cannot use such curves to prove that
$X$ is approximately elliptically connected  in the  cohomological sense.

We have however  the following result:
 \begin{lemm} $X$ is approximately elliptically connected (for any
 projective embedding $X\subset \mathbb{P}^N$) in the  cohomological sense.
  \end{lemm}
  {\bf Proof.} Indeed, recall that  $\tau:X
  \rightarrow A\times A\times \mathbb{P}^1$ is the blow-up map. Choose a neighborhood
  $\widetilde{U}$ of $\mathbb{P}(T_X)$ in $\mathbb{P}(T_{\mathbb{P}^N})$ and let $x\in X$.
  Let $\tau(x)=(y,t)$ with $y\in A\times A$ and $t\in \mathbb{P}^1$. We choose $x$ so that $t\not\in\{t_1,\ldots,t_4\}$, so that in particular
  $x$ does not belong to the exceptional divisor of $\tau$, and
  there is a copy $A\times A\times t\subset X$ passing through $x$. Then
  by theorem \ref{propexample}, $A\times A\times t$ is approximately elliptically connected in
   $\mathbb{P}^N$ in the infinitesimal and cohomological sense. Thus there is an elliptic curve $E_x$ passing through $x$, whose tangent lift
   is  contained
   in $\widetilde{U}$. Furthermore the class of these elliptic curves can be chosen to be independent
   of $x\in X$ and a convex combination of them generate an ample class on
   $A\times A\times t\subset X$.
   On the other hand, assume now that the line $y\times \mathbb{P}^1$ does not meet the  locus of $A\times A\times \mathbb{P}^1$  blown-up
   under  $\tau$. Then this line is
   a rational curve $C_x$ contained in $X$ and passing through $x$. In $\mathbb{P}^N$, we can smooth
   the curve $E_x\cup_x C_x$ and it is easy to see that the smoothed curve
   can be chosen to pass through $x$  and to stay
   $\widetilde{U}$-close to $X$. This proves
  the result since a convex combination of the classes $[C_x]+[E_x]$ is strongly positive.

\cqfd
\section{Stability results and further questions\label{sec2}}
We start with the following results concerning stability under \'{e}tale covers.
Here $P$ is any  smooth complex projective variety.
\begin{prop}\label{letale}
Assume $X\subset P$ is connected and  approximately rationally or elliptically connected  in the na\"{\i}ve,  resp. cohomological  sense.
Let $U\subset P$ be a neighborhood of $X$ which has the same homotopy type as $X$ (eg a tubular neighborhood). Then  any \'{e}tale connected proper cover
${X}'\rightarrow X$ is approximately rationally connected in the corresponding
neighborhood
$f:{U}'\rightarrow U$ of $U$, in the na\"{\i}ve,  resp. cohomological  sense.
\end{prop}
{\bf Proof.} We give the proof for the rational case, the elliptic case
working similarly, due to the fact that \'{e}tale covers of elliptic curves
are again elliptic curves.
Let us first consider approximate connectedness in the na\"{\i}ve sense.
Any small neighborhood $V_\epsilon({X}')$ of ${X}'$
in ${U}'$ contains a neighborhood of the form $f^{-1}({U_{\epsilon'}(X)})$, for some
$\epsilon'$. Let now ${x}',\,{y}'\in {X}'$ and let $x,\,y$
be their images in $X$. There is a smooth rational curve $C$ contained
in $U_{\epsilon'}(X)$, (we assume here ${\rm dim}\,U\geq3$ since the case where ${\rm dim}\,X=1$
is completely understood by Lemma \ref{lekoba},) and containing
$x$ and $y$. The inverse image of $C$ in $f^{-1}({U_{\epsilon'}(X)})$ is a finite  union of rational curves, and  one of them, say $C_{{x}'}$, passes through ${x}'$.
As $C_{{x}'}$ maps onto $C$, it contains one point ${y}''$
of ${X}'$ over $y$. In conclusion, the set of points
${y}''\in {X}'$ which are joined to ${x}'$
by a rational curve in
$V_\epsilon({X}')$ contains an open subset $W_{{x}'}\subset {X}'$ which maps onto $X$.
For any point ${z}\in W_{{x}'}$, the open set $W_{{z}}$ must be equal to
$W_{{x}'}$, since a point in $W_{{z}}$ is joined to ${x}'$ by a chain of two rational curves
passing through $z$, and this chain can be smoothed.
We may assume that the cover $f:{X}'\rightarrow X$ is Galois. Let $g\in {\rm Gal}\,({X}'/X)$. Then $gW_{{x}'}=W_{g{x}'}$, and by the above we conclude that
${X}'$ is the finite disjoint union of open sets of the form
$W_{{x}'}$. As ${X}'$ is connected, it follows that ${X}'=W_{{x}'}$.

For the approximate rational connected in the  cohomological sense, we have to add the following argument: recall
from Definition \ref{defi5} that we need to have for any tubular neighborhood $V'$ of ${X}'$, any neighborhood
$\widetilde{V}'$ of $\mathbb{P}(T_{X'})$ in $\mathbb{P}(T_{U'})$,
and any point ${x}'\in {X}'$ a finite number of  rational curves $C_{i,{x}'}$ passing through ${x}'$, contained in
$V'$, $\widetilde{V}'$-close to $X'$,  such that the class of the curve
$C_{i,x'}$ does not depend on $x'$ and the $(n-1,n-1)$-part of the class $\sum_in_i[C_{i,{x}'}]\in H_2(V',\mathbb{Z})=H_2({X}',\mathbb{Z})$ is strongly
positive for some $n_i>0$. Of course we may assume that $V'$ and $\widetilde{V}'$ are inverse images under $f$ of similar
neighborhoods $V$, $\widetilde{V}$ for $X\subset U$.
If we start from such data $C_{i,x}$ for $x\in X$ and for the neighborhoods $V$, $\widetilde{V}$, we observe now that the class in $H_2({X}',\mathbb{Z})$ of the unique
lift $C_{i,{x}'}$  of $C_{i,x}$ passing through ${x}'$  does not depend on ${x}'$, because it
can be chosen to vary
continuously with ${x}'$ and ${X}'$ is connected. In particular, we find that the sum
$\sum_{g\in {\rm Gal}\,({X}'/X)} g_*[C_{i,{x}'}]$ is equal to ${\rm card}\,G\,[C_{i,{x}'}]$ and it is the pull-back under $f$
of the class $[C_{i,x}]\in  H_2(X,\mathbb{Z})$. The fact that there exist such $C_{i,{x}'}$'s with $\sum_in_i[C_{i, x'}]^{n-1,n-1}\in H_2(V',\mathbb{R})=H_2({X}',\mathbb{R})$  strongly
positive is thus equivalent to the fact that there exist such $C_{i,{x}}$'s with $\sum_in_i[C_{i,{x}}]\in H_2(V,\mathbb{R})=H_2({X},\mathbb{R})$  strongly
positive.
\cqfd

The following consequence of Proposition  \ref{letale}  illustrates the power of the cohomological condition in
Definition \ref{defi6} :
\begin{coro} \label{coro24maifin} The varieties $X$ in Example \ref{ex23mai} are not approximately elliptically connected
in the cohomological sense in any projective embedding.
\end{coro}
{\bf Proof.} Recall that $X$ is a quotient of a product
$S\times C$ by a free involution $\iota$, where $g(C)\geq 2$.
If it was approximately elliptically connected
in the cohomological sense in some projective embedding, by Proposition \ref{letale}, the product
 $S\times C$ would be approximately rationally connected
in the cohomological sense in some  embedding. In particular, there would be
elliptic curves $E_i$ contained in a tubular neighborhood $U$ of $S\times C$, with the
property that some combination of the $(2,2)$-component of the classes
$[E_i]\in H_2(U,\mathbb{Z})=H_2(S\times C,\mathbb{Z})=H^4(S\times C,\mathbb{Z})$ is strongly positive.
But for any continuous map $\phi$ from an elliptic curve $E$ to a genus $\geq2$ curve $C$,
the induced map $\phi_*:H_2(E_i)\rightarrow H_2(C)$ is trivial. Thus the classes $pr_{2*}[E_i]$ vanish in $H_2(C,\mathbb{Z})$ and  for any line bundle $L$ of positive degree on $C$, $<pr_2^*c_1(L),[E_i]>=0$. Thus
$\sum_i n_i<pr_2^*c_1(L),[E_i]^{2,2}>=0$ for any $n_i$'s, which provides a contradiction.

\cqfd

Concerning the stability under morphism, we have the following easy result:

\begin{prop} \label{theosurj} Let $\phi: X\rightarrow Y$ be a surjective morphism, where $X$ and $Y$ are smooth
projective and ${\rm dim}\,Y>0$. If $X$ is  strongly approximately rationally (resp. elliptically) connected,  $Y$ is approximately rationally (resp. elliptically) connected in the na\"{\i}ve sense in any projective embedding.

\end{prop}
{\bf Proof.}   Let $j_Y: Y\hookrightarrow \mathbb{P}^N$ be a projective embedding.
Choose a projective  embedding $j_X$ of $X$ in some projective space $\mathbb{P}^M$,
and consider the embedding $j'_X=(j_X,j_Y\circ \phi):X\hookrightarrow P=\mathbb{P}^M\times \mathbb{P}^N$.
 By assumption  $j'_{X}(X)$
 is approximately rationally (resp. elliptically) connected in the na\"{\i}ve sense in
$P$. The morphism $pr_2:P\rightarrow \mathbb{P}^N$
 sends rational  curves (resp. chain of elliptic curves) passing through any two given points
 of $j'_X(X)$ and contained
in a sufficiently small neighborhood of $j'_X(X)$ to  rational  curves (resp. chain of elliptic curves)
passing through any two given points of $j_Y(Y)$
and contained in a given  neighborhood of ${Y}$ in $\mathbb{P}^N$.
\cqfd

\begin{coro} \label{coro3mai} A  fibration $X\rightarrow Y$ over a Kobayashi hyperbolic variety
$Y$  is not strongly
 approximately  elliptically  connected.
\end{coro}
{\bf Proof.} Indeed, if it was strongly approximately  elliptically  connected, the variety
$Y$  would be approximately  elliptically  connected, hence in particular
 not Kobayashi hyperbolic by Lemma \ref{lekoba}.  This gives a contradiction.
\cqfd

As we mentioned in the introduction,  Proposition \ref{theosurj} is
not true for na\"{\i}ve approximate rational or elliptic connectedness. This implies a negative answer to the following question:
\begin{question}\label{qversegre} Let $Z\subset \mathbb{P}^N$ be   the Segre embedding of $\mathbb{P}^k\times \mathbb{P}^l$ for some
 integers $k,\,l$. Fix a distance $d$ on $\mathbb{P}^N$. Is it true that for any $\epsilon>0$, there exists $\eta(\epsilon)>0$,
 such that $\lim_{\epsilon\rightarrow0}\eta(\epsilon) =0$, and that for any rational (resp. elliptic)
 curve $C$ contained in $U_\epsilon(Z)$, there is a rational (resp. elliptic) curve $C'\subset Z$
 such that $ d(C,C'):={\rm Sup}_{c\in C,c'\in C'} \{d(c,C'),d(c',C)\}\leq \eta$?
\end{question}
If we look at the proof of Lemma \ref{exa10mai}, i), a counterexample is obtained by constructing in a small neighborhood
of the union of
a large number of lines $l_i=x_i\times \mathbb{P}^1,\,i=1,\ldots,\,M$ contained in  $Z$, with $d(l_i,l_{i+1})<\epsilon$,  a chain $\cup_il_{i,\epsilon}$ of rational  curves in $\mathbb{P}^N$ obtained by deforming the $l_i$ in such a way
that $l_{i,\epsilon}$ meets $l_{i+1,\epsilon}$, and then by smoothing the resulting chain to a
rational curve $C$. If the diameter of the set $\{x_i,\,i=1,\ldots,\,M\}$ is large, and the points $x_i$ are taken
in a Kobayashi hyperbolic subvariety $Y\subset \mathbb{P}^k$, the distance $d(C,C')$
 between $C$ and any  elliptic or rational  curve  $C'$ in $Z$ is bounded below by a positive constant.

\subsection{Further questions and remarks}

The first obvious question is the following:
\begin{question} Let $X\subset \mathbb{P}^N$ be  approximately elliptically  connected (in the strong or cohomological sense). Is the Kobayashi
pseudo-distance of $X$ trivial?
\end{question}
As our motivation was to understand the class of varieties with trivial Kobayashi pseudodistance, which includes conjecturally Calabi-Yau manifolds (cf. \cite{koba}), it is also natural to ask
the following:
\begin{question} Let $X$ be a Calabi-Yau manifold. Is $X$ approximately rationally or elliptically
connected (in any of the senses introduced in this paper)?
\end{question}
Another important question   is the following:
\begin{question} Is the property of cohomological approximate rational or  elliptic connectedness independent of the choice of projective embedding?
\end{question}

The following  question is related to the work of Graber, Harris and Starr \cite{GHS}:
\begin{question}\label{questionfiberbase} Let $\phi:X\rightarrow Y$ be a surjective morphism. Assume that the fibers
of $\phi$ are rationally connected (see \cite{komimo}) and that $Y$ is
approximately rationally (resp. elliptically)  connected in the strong or cohomological sense. Then is  $X$
approximately rationally (resp. elliptically)  connected in the same sense?

\end{question}
Let us give one result in this direction: let $P=\mathbb{P}^{n}$ and let
$Q=\mathbb{P}(H^0(P,\mathcal{O}_P({d})))$, where
${d}\leq n$ and $n\geq 2$  so that the general
 hypersurface of degree  ${d}$ in $P$  is Fano of dimension $\geq1$. There is a universal subvariety $Z\subset Q\times P$, defined by the  tautological
equation $F_{{d}}\in H^0(\mathcal{O}_{Q\times P}(1,{d}))$.  Via the second projection, $Z$ is a fibration in projective spaces over $P$.
\begin{prop}\label{fibreci} If $Y\subset Q$ is rationally (resp. elliptically) connected in the cohomological sense, and
 $X:=Y\times_ QZ\rightarrow Y$ is smooth of the expected dimension (hence the generic fiber of $X\rightarrow Y$ is a smooth
 hypersurface in $P$),
 $X$ is approximately rationally (resp. elliptically) connected in the cohomological sense in $Z$, hence in  $Q\times P$.
\end{prop}
{\bf Proof.} Let $V\subset Z$ be a tubular neighborhood of $X$ and $\widetilde{V}\subset \mathbb{P}(T_Z)$ be a neighborhood of
$\mathbb{P}(T_X)$. There are neighborhoods $U\subset Q$ of $Y$ and $\tilde{U}\subset \mathbb{P}(T_Q)$ of $\mathbb{P}(T_Y)$
such that $V$ contains $\pi^{-1}(U)$ and $\widetilde{V} $ is contained in $\pi_*^{-1}(\tilde{U})$ where $\pi:={pr_{1\mid Z}}:Z\rightarrow Q$.
As $Y\subset Q$ is  approximately rationally, resp. elliptically, connected in the cohomological sense, there is a  curve $E$ which is rational (resp.  elliptic),  contained in $U$ and passing through any  point
$y$ of $Y$. Furthermore, $E$ can be chosen to be $\widetilde{U}$-close to $Y$, of class independent
of $y$, and finally a convex combination of the $(n-1,n-1)$-part of these classes contains a strongly
positive class in $Y$, where $n={\rm dim}\,Y$.
Moving $E$ if needed, we may assume
that  $Z_E$ is smooth with irreducible fibers and $Z_E\rightarrow E$ is a smooth Fano complete intersection over the generic point of  $E$.
By the Tsen-Lang theorem, the family $Z_{E}\rightarrow E$ has a section. Results of \cite{komimo}
even show that such a section $\widetilde{E}$ can be chosen to have an arbitrarily positive class in $Z_E$.
These sections produce elliptic curves $\widetilde{E}\subset V$ which are then $\widetilde{V}$-close to $X$,
and  pass through the general point of $X$. Finally, under our assumptions, (and because we may assume that ${\rm dim}\,X\geq3$, otherwise the result is obvious) the Lefschetz hyperplane section theorem says
that $H^2(X,\mathbb{Z})=H^2(Y,\mathbb{Z})\oplus H^2(P,\mathbb{Z})$. It is then immediate
to conclude that if the curves $\widetilde{E}$ have an ample class in $Z_E$ and
a convex combination of the $(n-1,n-1)$-components of the pushforward of their classes in
$H_2(U,\mathbb{Z})=H_2(Y,\mathbb{Z})=H^{2n-2}(Y,\mathbb{Z})$ contains a strongly positive class, then
a convex combination of the $(m-1,m-1)$-components of  their classes in
$H_2(V,\mathbb{Z})=H_2(X,\mathbb{Z})=H^{2m-2}(X,\mathbb{Z})$ contains a strongly positive class.

\cqfd
\begin{rema} {\rm The analogous result, if one only assumes that the fibers of $X\rightarrow Y$ are  approximately elliptically or even rationally   connected in the strong sense,  is not true.
Indeed, consider the example \ref{ex23mai} where  $X=(C\times S)/\iota$, and $Y=\mathbb{P}^1$,
where $C$ is a curve of genus $\geq 2$ with hyperelliptic involution $i$, $S$ is a $K3$
surface which is the universal cover $S\rightarrow T$ of an Enriques surface.
 The morphism $\phi:X\rightarrow Y$ is induced by passing to the quotient from the projection $p_2:C\times S\rightarrow C$ using the isomorphisms
  $X\cong (C\times S)/\iota,\,C/i\cong \mathbb{P}^1$
  and equivariance of  $p_2$. The fibers of $\phi$ are isomorphic
to $S$ or to $T$, hence are strongly approximately rationally connected. However  $X$ is not
strongly approximately rationally  or elliptically connected by Corollary \ref{coro24maifin}.

}
\end{rema}
To finish, let us conclude with the following questions:
\begin{question} \label{qratpi1}(Campana)  Assume $X$ is  approximately rationally connected (in the adequate sense). Is $\pi_1(X)$ finite?
\end{question}

The following similar  question is very much related to the results of \cite{campana}.

\begin{question} \label{qellpi1} (Campana)  Assume $X$ is strongly  approximately elliptically connected (in the adequate sense). Is $\pi_1(X)$ virtually abelian?
\end{question}

The two questions (for cohomological approximate connectedness) are related as follows.
\begin{prop} Suppose Question \ref{qellpi1} has
a positive answer for  cohomological approximate connectedness, then Question \ref{qratpi1} also has a positive answer for cohomological approximate connectedness.
\end{prop}
{\bf Proof.} Let $X\subset \mathbb{P}^N$ be approximately rationally connected in the cohomological sense.
 We know, assuming Question \ref{qellpi1} has
a positive answer, that $\pi_1(X)$
is virtually abelian. Passing to an \'{e}tale cover of $X$, and using
 Lemma \ref{letale}, we may assume that $X$ is  approximately rationally connected in the
  cohomological sense in an adequate variety $U$, and furthermore has torsion free abelian $\pi_1$. We  want  to prove that
  $\pi_1({X})$ is  trivial.   Equivalently, if $a_X:X\rightarrow Alb\,X$ is the Albanese map,  letting $Y:=a_X(X)\subset Alb\,X$,
 one wants to prove that $Y$ is a point. Assume the contrary. Then choosing an ample line bundle on
 $Alb\,X$, its pull-back $a_X^*L$ to $X$ is a semi-positive line bundle which is not numerically trivial.
 Consider now rational curves $C$ in a tubular neighborhood $U$ of $X$ in projective space.
 Then their class $[C]\in H_2(U,\mathbb{Z})\cong H_2(X,\mathbb{Z})$
 factors through $\pi_2(X)$, hence vanishes in $H_2(Alb\,X)$ under the map $a_{X*}$.
 Hence we conclude that $<[C],\phi^*c_1(a_X^*L)>=0$.

 This contradicts the fact that $X$ is approximately rationally connected in the cohomological sense,
 because the later implies in particular the existence of rational curves $C_i$ in any small neighborhood of $X$ in $U$, with the property that the class
 $\sum_i n_i[C_i]^{n-1,n-1}\in H^{n-1,n-1}(X)$ is strongly positive, so that $$\sum_in_i<[C_i],\phi^*c_1(a_X^*L)>\not=0.$$
 \cqfd

\end{document}